\documentclass{amsart}

\allowdisplaybreaks
\usepackage{mathpazo}
\usepackage{amsmath,amsfonts,amssymb}
\usepackage{amsrefs}
\usepackage{amsthm}
\usepackage{latexsym,amsmath,amssymb,amsfonts}
\usepackage{rotating}
\usepackage{mathrsfs}
\usepackage{xypic} \xyoption{all}
\usepackage{amscd}
\usepackage{hyperref} 
\usepackage{euscript}
\usepackage{hhline}
\usepackage{graphicx,epstopdf}
\usepackage{epsfig}
\usepackage{xcolor}
\usepackage[all,color]{xy}

\newlength{\fighskip} \fighskip=2pt
\newlength{\figvskip} \figvskip=3pt

\usepackage{hyperref}

\advance\textheight by \topskip
\numberwithin{equation}{section}

\newcommand{\op}{\operatorname}
\newcommand{\C}{\mathbb{C}}

\newcommand{\R}{\mathbb{R}}

\newcommand{\Z}{\mathbb{Z}}
\newcommand{\Etau}{{\text{E}_\tau}}
\newcommand{\E}{{\mathcal E}}

\renewcommand{\H}{\mathbb{H}}

\newcommand{\g}{\mathfrak{g}}
\newcommand{\im}{\op{im}}

\providecommand{\abs}[1]{\left\lvert#1\right\rvert}

\newcommand{\abracket}[1]{\left\langle#1\right\rangle}
\newcommand{\bbracket}[1]{\left[#1\right]}
\newcommand{\fbracket}[1]{\left\{#1\right\}}
\newcommand{\bracket}[1]{\left(#1\right)}

\newcommand{\mc}{\mathcal}
\newcommand{\cinfty}{C^{\infty}}
\newcommand{\pa}{\partial}

\renewcommand{\dbar}{\bar\pa}
\newcommand{\OO}{{\mathcal O}}

\newcommand{\BV}{Batalin-Vilkovisky}

\newcommand{\iso}{\cong}

\DeclareMathOperator{\Sym}{Sym}
\DeclareMathOperator{\Hom}{Hom}

\DeclareMathOperator{\Tr}{Tr}

\DeclareMathOperator{\Obs}{Obs}

\DeclareMathOperator{\PV}{PV}

\DeclareMathOperator{\Jac}{Jac}
\newcommand{\A}{\mathcal A}

\theoremstyle{plain}
\newtheorem{thm}{Theorem}[section]
\newtheorem{thm-defn}{Theorem/Definition}[section]
\newtheorem{lem}[thm]{Lemma}
\newtheorem{lem-defn}[thm]{Lemma/Definition}

\newtheorem{cor}[thm]{Corollary}

\theoremstyle{definition}
\newtheorem{defn}[thm]{Definition}

\newtheorem{eg}[thm]{Example}

\theoremstyle{remark}
\newtheorem{rmk}[thm]{Remark}

\allowdisplaybreaks[4]  

\begin{document}

 \title{Effective Batalin-Vilkovisky quantization and geometric applications}
  \author{Si Li}
  \date{}

  \maketitle


\begin{abstract} We explain the effective renormalization method of quantum field theory in the \BV\  formalism and illustrate its mathematical applications by three geometric examples:  (1) Topological quantum mechanics and algebraic index, (2) Elliptic curve and higher genus mirror symmetry, (3) Calabi-Yau geometry and integrable hierarchy. This note is an expansion of author's talk at ICCM 2016. 
\end{abstract}

  \tableofcontents


\section{Introduction}
Quantum field theory has been a rich source of mathematical structures. Typically,  quantum field theory deals with infinite dimensional geometry, a fact that lies behind many of its nontrivial consequences. A famous example is the mysterious mirror symmetry conjecture between symplectic and complex geometries. It has a clean physics interpretation as a version of infinite dimensional Fourier transform. Still, we do not have a good mathematical framework to understand why mirror symmetry holds although many nontrivial statements have been checked.  Nevertheless, asymptotic analysis can always be performed with the help of the celebrated idea of renormalization and effective field theory. The purpose of this article is to explain some aspects of effective renormalization method of quantum field theory and show how it leads to powerful applications to geometric problems. 

A physics system is usually described by a map
$$
  S: \mc E \to \R. 
$$
Here $\mc E$ is called the \emph{space of fields} and $S$ is the \emph{action functional}.  Here are some typical examples. 
\begin{itemize}
\item $\mc E=C^\infty(X)$. \emph{Scalar field theory}. 
\item $\mc E=\{\text{connections on a vector bundle}\ V\to X\}$.  \emph{Gauge theory}. 
\item $\mc E=\text{Map}(\Sigma, X)$. \emph{$\sigma$-model}. 
\item $\mc E=\{\text{metrics on}\ X\}$. \emph{Gravity}. 
\end{itemize}

In all examples above, $\E$ is an infinite dimensional space carrying a delicate topology. Classical physics is decribed by the critical locus
$$
   \text{Crit}(S)=\{\delta S=0\}
$$
which usually leads to interesting differential equations (\emph{equations of motion}). One convenient way to approach quantum physics is by Feynman's {``path integral"}
$$
\abracket{\mc O}:=   \int_{\mc E}  \mc O e^{S/\hbar}.
$$   
$\mc O$ is called an \emph{observable}. $\abracket{\mc O}$ is the \emph{correlation function}. Again, we encounter the mathematical challenge for such $\infty$-dim integral. On the other hand, $\hbar$-asymptotic analysis has a solid foundation usually called \emph{perturbative field theory}. Although modern physics is deeply rooted in understanding nonperturbative aspects of quantum field theory, the $\hbar$-asymptotic analysis gives very useful physics information as well as new mathematical phenomenons. We illustrate this by three geometric applications. 

\begin{enumerate}
\item \textbf{Topological quantum mechanics and index theorem}. Quantum field theory in one dimension is quantum mechanics. We first explore the well-studied example of topological quantum mechanics and explain the result of \cite{GLL} on its connection with Fedosov's deformation quantization \cite{Fedosov-DQ}:
$$
\xymatrix{
\boxed{\text{low energy effective theories of topological quantum mechanics }}\ar@{=>}[d]\\
\boxed{\text{Fedosov's abelian connections on the Weyl bundle}}
}.
$$
This fact together with a version of Atiyah-Bott localization computation leads to a geometric proof of the \emph{algebraic index theorem}, which was first formulated by Fedosov \cite{Fedbook} and Nest-Tsygan \cite{Nest-Tsygan} as the algebraic analogue of Atiyah-Singer index theorem. This result is analogous to  the physics interpretation of Atiyah-Singer index theorem \cite{Alvarez, FW,Witten, Witten-index} via semi-classical analysis of supesymmetric quantum mechanics, as well as the symbol calculus of pseudodifferential operators \cite{G-index, Nest-Tsygan2,Fedesov-PDO}.
 
\item \textbf{Elliptic curve and higher genus mirror symmetry}. We consider chiral deformations of free CFT's in the \BV\ (BV) formalism in two dimension \cite{Si-vertex}. We explain that 
$$
\xymatrix{
\boxed{\text{BV quantization of 2d chiral theory}}\ar@{=>}[d]\\
\boxed{\text{Integrable evolution of vertex algebra}}
}
$$
where the evolution is parametrized by the relevant background symmetry. This can be viewed as the 2d analogue of the above 1d quantum mechanical result. As an application, we show how this construction leads to an exact solution \cite{Si-vertex} of quantum B-model (BCOV theory \cite{Si-Kevin, BCOV}) in complex one dimension that solves the higher genus mirror symmetry conjecture on elliptic curves \cite{Si-vertex, L-elliptic}.

\item \textbf{Calabi-Yau geometry and integrable hierarchy}. Topological string on Calabi-Yau geometry usually leads to integrable hierarchies. We propose a geometric interpretation of this mysterious phenomenon in terms of B-twisted topological string field theory (BCOV theory \cite{Si-Kevin, BCOV}). Let $X$ be a Calabi-Yau geometry (either a compact Calabi-Yau or a Landau-Ginzburg model). We consider $X\times \C$ which is again Calabi-Yau but with dimension increased by one. Compactifying along $X$, 
$$
\xymatrix{
\boxed{\text{BCOV theory on}\ X\times \C}\ar@{=>}[d]\\
\boxed{\text{2d chiral theory on}\ \C}
}
$$
we get an effective 2d chiral theory on $\C$. We explain the proposal in \cite{HSY} that the BV master equation for the effective 2d chiral theory induced from that on $X\times\C$ leads to integrable hierarchies.
$$
   \boxed{\text{BV master equation}} \Longrightarrow \boxed{\text{Integrable hierarchy}}
$$ 
\end{enumerate}

In all examples above, we work with Batalin-Vilkovisky (BV or BRST-BV) quantization scheme, which arises in physics as a general framework to quantize theories with gauge symmetries \cite{BV}. In this article we work with Costello's formalism\cite{Kevin-book} of homotopic effective renormalization which is a convenient framework for our geometric problems. We introduce the basic philosophy and notions of BV quantization in the asymptotic analysis of quantum gauge theories, and then explain with more details on the above examples. 

\noindent \textbf{Conventions}

In this article we work with the category of $\Z$-graded vector spaces over a field of characteristic $0$. 
\begin{itemize}
\item Let $V=\bigoplus\limits_{m\in \Z} V_m$ be a $\Z$-graded $k$-vector space. Given $a\in V_m$, we let $\bar a=m$ be its degree.
\item $V[n]$ denotes the degree shifting such that $V[n]_m=V_{n+m}$.
\item All tensors are in the category of graded vector spaces. Koszul sign rule is assumed. $\Sym^m(V)$ and $\wedge^m(V)$ denote the graded symmetric product and graded skew-symmetric product respectively.
\item $V[\hbar]$, $V[[\hbar]]$ and $V((\hbar))$ denote polynomial series, formal power series and Laurent series respectively in a variable $\hbar$ valued in $V$.

\item Given $A, B\in \Hom(V, V)$, the commutator $[A, B]\in \Hom(V, V)$ always means the graded commutator: $[A, B]=AB-(-1)^{\bar A \bar B}BA$. 

\end{itemize}

\section{\BV\ quantization}

\subsection{A prototype of BV formalism}
Calculus is nicely packaged into the framework of de Rham theory: let $X$ be an oriented manifold, $C\subset X$ be a submanifold, $\alpha$ be a smooth differential form, then there is a natural integration pairing between $\alpha$ and $C$
$$
  \int_{C}\alpha. 
$$
If $\alpha$ is closed, i.e. $d\alpha=0$, and $C$ has no boundary, then Stokes' theorem tells that the above integration is invariant under continuous deformations of $C$. 

The BV approach to the above problem is to consider the super manifold $T^*X[-1]$, where we have graded the cotangent fiber direction to be odd of degree $1$.  It carries a natural $(-1)$-shifted symplectic structure (which we denote by $\omega$), and we can identify its structure sheaf via polyvector fields 
$$
  \OO_{T^*X[-1]}=\PV^{-\bullet}(X), \quad \text{where}\ \PV^{-k}(X):=\Gamma(X, \wedge^k T_X). 
$$
The induced Poisson bracket on $\OO_{T^*X[-1]}$ coincides with the Schouten-Nijenhuis bracket.

Let us fix a chosen volume form $\Omega$ on $X$. It allows us to identify via contracting with $\Omega$
$$
  \PV^{-\bullet}(X)\stackrel{\lrcorner \Omega}{\longleftrightarrow} \Omega^{n-\bullet}(X), \quad n=\dim(X).
$$
Under this identification, the de Rham differential $d$ gives rise to a differential 
$$
  \Delta: \PV^{-\bullet}(X)\to \PV^{1-\bullet}(X).
$$
$\Delta$ will be called the BV operator, which is the divergence operator with respect to the volume form $\Omega$. The triple $(T^*X[-1], \omega, \Delta)$ is an example of BV manifold \cite{Schwarz}. In general, a BV manifold $M$ carries a natural integration pairing between its functions and lagrangian super submanifolds \cite{Schwarz}:
$$
    \int_{\mc L}: \OO_M\to \R. 
$$
If $\mc L$ has no boundary, and $\mu$ is $\Delta$-closed, then $\int_{\mc L}\mu$ is invariant under continuous deformations of $\mc L$ and a version of Stokes' theorem holds. This is the BV analogue of de Rham theory.

In our case 
$$
  \int_{\mc L}:  \OO_{T^*X[-1]}\to \R,
$$
if $\mc L=N^*_{C/X}[-1]\subset T^*X[-1]$ is the (shifted) conormal bundle of a submanifold $C$ inside $X$, then 
$$
  \int_{\mc L}\mu= \int_C \mu\lrcorner \Omega. 
$$
In particular, $\int_{\mc L}\Delta\mu=\int_{C}d(\mu \lrcorner \Omega)=0$ if $\pa C=\emptyset$. 

In our geometric situation above, the de Rham theory and BV formalism are completely equivalent. There is one subtle difference when we move on to discuss infinite dimensional geometry of quantum field theories.  In de Rham theory, the integration map 
$$
   \int_X:  \Omega^{\bullet}(X)\to \R
$$
can be identified with the cohomology 
$$
  \int_X\iso H^n:  \Omega^{\bullet}(X)\to  H^{n}(X)\iso \R. 
$$ 
Here $n$ is the dimension of $X$. In the BV approach, the BV integration 
$$
  \int_X: \PV^{-\bullet}\to \R
$$  
can be identified with $H^0$. When $n=\infty$ as in quantum field theory, $H^0$ behaves better than $H^n$. The problem of integration is transferred to construct $\Delta$, which has a convenient formulation at least in perturbative quantum field theory. This is one of the folklore advantage of working with BV. 
\subsection{BV algebra and master equation}
Let us first introduce some algebraic structures that are relevant to the above discussion. 

\begin{defn}\label{defn-BV} A differetial \BV\ (BV) algebra is a triple $(\A, Q, \Delta)$ where
\begin{itemize}
\item $\A$ is a $\Z$-graded commutative associative unital algebra,
\item $\Delta: \A \to \A$ is a second-order operator of degree $1$ such that $\Delta^2=0$,
\item $Q: \A\to \A$ is a derivation of degree $1$ such that $Q^2=0$ and $[Q, \Delta]=0$. 
\end{itemize}
\end{defn}

Here $\Delta$ is called the BV operator. $\Delta$ being ``second-order" means the following: let us define the \emph{BV bracket} $\fbracket{-,-}$ as the failure of $\Delta$ to be a derivation
$$
    \fbracket{a,b}:=\Delta(ab)-(\Delta a)b- (-1)^{\bar a}a \Delta b. 
$$
Then $\fbracket{-,-}$ defines a Lie bracket of degree $1$ (Gerstenhaber algebra) such that $\Delta$ is compatible with $\{-,-\}$ via a graded version of Leibniz rule. In our discussion above, the triple is $(\PV^{-\bullet}(X), 0, \Delta)$.

\begin{defn} Let $(\A, Q, \Delta)$ be a differential BV algebra. An element $
I=I_0+I_1\hbar+\cdots \in \mc A[[\hbar]]$ of degree 0 is said to satisfy \emph{quantum master equation (QME)} if 
$$
QME: \quad \boxed{ (Q+\hbar \Delta) e^{I/\hbar}=0}. 
$$
This is equivalent to 
$$
\boxed{QI+\hbar \Delta I+{1\over 2}\{I,I\}=0}. 
$$
The leading $\hbar$-order $I_0$ satisfies 
$$
 CME: \quad \boxed{QI_0+{1\over 2}\{I_0, I_0\}=0}
$$
which is called the \emph{classical master equation (CME)}. 
\end{defn}

\begin{defn}A solution $I$ of quantum master equation gives rise to a differential 
$$
Q+\hbar\Delta+\fbracket{I,-}: \A[[\hbar]]\to \A[[\hbar]].
$$
The cohomology 
$$
\Obs^q=\boxed{H^{\bullet}(\A[[\hbar]], Q+\hbar \Delta+\{I,-\})}
$$ 
is called the \emph{quantum observables}. Similarly, its classical limit 
$$
\Obs^c=\boxed{H^{\bullet}(\A, Q+\{I_0,-\})}
$$ 
is called the \emph{classical observables}. 
\end{defn}

Here we give some mathematical and physical explanations of the above definitions. 
\begin{enumerate}
 
 \item {Classical master equation}: $
   Q+\{I_0,-\}
$
squares zero and describes the infinitesimal gauge transformations (BRST). \mbox{Mathematically, this defines a homotopic Lie-algebra ($L_\infty$-algebra)} $\g$ if $\mc A=\prod\limits_{n\geq 0}\Sym^n(\g[1]^*)$.
\item {Quantum master equation}:   the quantum gauge consistency condition for the measure $e^{I/\hbar}$.

\item A {quantum observable} is represented by an element $\mc O\in \mc A[[\hbar]]$ such that
$$
   (Q+\hbar \Delta)(\mc O e^{I/\hbar})=0. 
$$
This is the analogue of a closed differential form in the BV setting. 
\item {Quantum correlation function} is given by a BV integral
$$
  \abracket{\mc O}= \int_{\mc L} \mc O e^{I/\hbar}. 
$$
$\mc L$ is a lagrangian super subspace of $\mc E$ related to a \emph{gauge fixing condition}. In particular, for a quantum observable $\mc O$, $\abracket{\mc O}$ is invariant under continuous deformations of the gauge fixing condition $\mc L$ and independent of the choice of representatives. This is the original motivation for the BV description of quantum gauge consistency. 
\end{enumerate}

\begin{eg}[Singularity theory] We compare  BV formalism in the finite dimensional case with classical singularity theory. Let $f: \C^n\to \C$  be a polynomial in $n$ variables with an isolated critical point at the origin. Let $\{z^i\}_{i=1}^n$ be holomorphic coordinates on $\C^n$. Let $\{\theta_i\}_{i=1}^n$ be Grassmann variables of degree $-1$
$$
  \theta_i \theta_j=-\theta_j\theta_i, \quad \forall 1\leq i, j\leq n. 
$$
Let $\C\{z\}$ be the germ of holomorphic functions at $0$. We consider the triple $(\A=\C\{z^i\}[\theta_i], Q=0, \Delta)$ where
$$
  \Delta=\sum_{i=1}^n {\pa\over \pa \theta_i}{\pa\over \pa z^i}.
$$
It is easy to check that $\Delta^2=0$, and its expression explains the name ``second order". 
\begin{itemize}
\item $f=f(z^i)$ gives a solution of QME: $\Delta e^{f/\hbar}=0$ since $f$ does not depend on $\theta_i$'s. 
\item The space of classical observables is given by the Jacobian algebra or Milnor ring
$$
   \Obs^c=H^{\bullet}(\A, \{f,-\})=\Jac(f):= \C\{z^i\}/(\pa_i f).
$$
\item The space of quantum observables is isomorphic to a formal completion of the Brieskorn lattice  \cite{Saito-residue}
$$
  \Obs^q= H^{\bullet}(\A[[\hbar]], \hbar\Delta+\{f,-\})\stackrel{d^n z}{\iso}  {\Omega^n(\C^n)_0[[\hbar]]\over (\hbar d+ df\wedge )\Omega^{n-1}(\C^n)_0[[\hbar]]}
$$
Here $d^n z=dz^1\wedge \cdots \wedge dz^n$, $\Omega^k(\C^n)_0$ is the germ of holomorphic $k$-forms at the origin $0$. 
\item The BV integration models the oscillatory integral
$$
     \abracket{\mc O}=\int_{\mc L}  d^nz \mc O e^{f/\hbar} 
$$
where $\mc L$ is a Lefschetz thimble. 
\end{itemize}
The Brieskorn lattice (or the quantum observables) plays the important role of Hodge filtration in singularity theory (see \cite{Kulikov, AGV} for an exposition). The above finite dimensional model is the effective theory of topological Landau-Ginzburg B-model. In general, a quantum field theory can be viewed as a version of infinite dimensional singularity theory, and the $\hbar$-filtration is related to the Hodge filtration. 
\end{eg}

\subsection{(-1)-shifted symplectic geometry}
 A general class of differential BV algebras arises from odd symplectic geometry.  Let us start with a finite dimensional toy model. Let $(V,Q)$ be a finite dimensional dg vector space. The differential $Q: V\to V$ induces a differential on various tensors of $V, V^*$, still denoted by $Q$. Let 
$$
   \omega \in \wedge^2 V^*, \quad Q(\omega)=0, 
$$
be a $Q$-compatible symplectic structure such that $\deg(\omega)=-1$. Then $\omega$ induces an isomorphism
$$
     V^* \simeq V[1].
$$
Let $K=\omega^{-1}\in \Sym^2(V)$ be the Poisson kernel of degree $1$ under the identification via $\omega$
\[     \wedge^2 V^* \simeq\  \Sym^2(V)[2],
\]
where we have used the canonical isomorphism $\wedge^2(V[1])\simeq \Sym^2(V)[2]$. Let us define the formal power series ring
$$
 \OO(V):=\prod_{n\geq 0}\Sym^n(V^*).
$$
Then $(\OO(V), Q)$ is a graded-commutative dga. Let $\Delta_K$ denote the second order operator
$$
   \Delta_K: \Sym^n(V^*)\to \Sym^{n-2}(V^*)
$$
by contracting with the kernel $K\in \Sym^2(V)$ using the pairing between $V$ and $V^*$. It is straight-forward to see that 
$(\OO(V), Q, \Delta_K)$ defines a differential BV algebra. 

\begin{rmk}\label{rmk-BV-degenerate}
The relevant data defining the BV operator is $K$ instead of $\omega$. In fact,  for any $Q$-closed element $K\in \Sym^2(V)$ of degree $1$, the triple $(\OO(V), Q, \Delta_K)$ defines a differential BV algebra. In particular, the above construction is generalized to the case when $K$ is degenerate, i.e. the Poisson case.  
\end{rmk}

Let $P$ be a degree $0$ element of $\Sym^2(V)$, and consider the new BV kernel 
$$
  K_P:= K+Q(P). 
$$
We say that $K_P$ is homologous to $K$. We can form a new differential BV algebra $(\OO(V), Q, \Delta_{K_P})$. Let $\pa_P$ denote the second order operator on $\OO(V)$ by contracting with $P$, which is defined similarly to $\Delta_K$ associated to $K$. The following simple algebraic lemma is key to Costello's homotopic renormalization theory \cite{Kevin-book}. 

\begin{lem}\label{lem-HRG}  The following equation holds formally as operators on $\OO(V)[[\hbar]]$
$$
     \bracket{Q+\hbar \Delta_{K_{P}}}e^{\hbar \pa_{P}}=e^{\hbar \pa_{P}}\bracket{Q+\hbar \Delta_{K}},
$$
i.e., the following diagram commutes
$$
\xymatrix{
    \OO(V)[[\hbar]] \ar[rr]^{Q+\hbar \Delta_{K}} \ar[d]_{\exp\bracket{\hbar \pa_{P}}} && \OO(V)[[\hbar]] \ar[d]^{\exp\bracket{\hbar \pa_{P}}}\\
   \OO(V)[[\hbar]] \ar[rr]_{Q+\hbar \Delta_{K_{P}}} && \OO(V)[[\hbar]]
}
$$
\end{lem}
This lemma is a simple consequence of  the observation that 
$
  \bbracket{Q, \pa_{P}}=\Delta_{K}-\Delta_{K_{P}}. 
$

\begin{defn} Let 
$$
\OO^+(V):= {\Sym^{\geq 3}(V^*))}\oplus \hbar \OO(V)[[\hbar]]
$$
denote the subspace of $\OO(V)[[\hbar]]$ consisting of functions that are at least cubic modulo $\hbar$. 
\end{defn}

\begin{cor}\label{cor-HRG} Let $I\in \OO^+(V)$ be a solution of quantum master equation in $(\OO(V)[[\hbar]], Q, \Delta_K)$. Let $I_P$ be formally defined by the equation
$$
  e^{I_P/\hbar}= e^{\hbar \pa_P} e^{I/\hbar}. 
$$ 
Then $I_P$ is a well-defined element of $\OO^+(V)$ that solves the quantum master equation in $(\OO(V)[[\hbar]], Q, \Delta_{K_P})$.
\end{cor}
The real content of the above formal definition of $I_P$ is that 
$$
  I_P= \sum_{\Gamma: \text{connected}} \frac{\hbar^{g(\Gamma)}}{|\text{Aut}(\Gamma)|} W_{\Gamma}(P, I) 
$$
where the summation is over all connected Feynman graphs with $P$ being the propagator and $I$ being the vertex (see \cite{graph} for an exposition on Feynman graph techniques). Here $g(\Gamma)$ denotes the genus of $\Gamma$, and $|\text{Aut}(\Gamma)|$ denotes the cardinality of the automorphism group of $\Gamma$. $I$ being at least cubic  implies that the above infinite graph sum leads to a well-defined element of $\OO^+(V)$. The statement of $I_P$ being a solution of the quantum master equation is a direct consequence of Lemma \ref{lem-HRG}.

\begin{defn} Given $P\in \Sym^2(V)$ of degree $0$,   we define the \emph{homotopic renormalization group flow} (HRG) operator $W(P, -)$ by
\begin{align*}
  W(P, -): \OO^+(V)\to \OO^+(V), \quad 
           \boxed{e^{W(P,I)/\hbar}:=e^{\hbar \pa_P} e^{I/\hbar}}. 
\end{align*}
\end{defn}

Lemma \ref{lem-HRG} says that HRG links QME for differential BV algebras whose BV kernels are homologous. 

\subsection{Effective BV quantization}
Quantum field theories always have infinite dimensional space of fields playing the role of $V$ in the previous subsection. Typically, a classical field theory in the BV formalism consists of $(\E, Q, \omega,  S_0)$ where
\begin{enumerate}
\item fields $\mc E= \Gamma(X, E)$ for a graded bundle $E$ on a manifold $X$. 
\item $Q: \mc E\to \mc E$  is a differential operator on $E$ that makes $(\E, Q)$ into an elliptic complex.
\item $\omega$: local symplectic pairing of degree $(-1)$
$$
  \omega(\alpha, \beta)=\int_X \abracket{\alpha, \beta}, \quad \forall \alpha,\beta\in \mc E. 
$$
where $(-,-)$ is a degree $-1$ skew-symmetric pairing on $\E$ valued in the density line bundle on $X$. 

\item A classical action $
  S_0=\omega(Q(-),-)+I_0
$ satisfying the classical master equation 
$$
  \{S_0,S_0\}=0, \quad \text{equivalently}\ QI_0+{1\over 2}\{I_0,I_0\}=0. 
$$ 
Here the quadratic term $\omega(Q(-),-)$ is called the \emph{free part}, while $I_0$ is called the \emph{interaction}. 
\end{enumerate}

\begin{eg}[Chern-Simons theory] Let $X$ be a 3-manifold. $\g$ be a Lie algebra with a Killing pairing $\Tr$. 
\begin{itemize}
\item $\mc E=\Omega^\bullet(X, \g)[1]$, and $Q=d$ is the de Rham differential. 
\item 
$
  \omega(\alpha, \beta)=\int_X \Tr (\alpha\wedge \beta), \quad \alpha, \beta\in \mc E
$.
\item  Let $[-,-]$ be the bracket on $\mc E$ induced from $\g$. The {Chern-Simons functional} in the BV formalism is
$$
 CS(A)={1\over 2}\int_{X} \Tr(A\wedge dA)+{1\over 6}\int_X \Tr(A\wedge [A,A]), \quad A\in \mc E
$$ 
\item Classical master equation $\Longleftrightarrow \fbracket{\Omega^\bullet(X,\g), d, [-,-]}$ being a {differential graded Lie algebra}. 
\end{itemize}

\end{eg}

The space of fields $\E$ carries a natural topology, and $V^*$ is replaced by the continuous dual of $\E$
$$
  \E^*=\Hom(\E, \R),
$$
i.e., distributions. Then formal functions $\OO(\E)$ are defined in terms of the completed tensor product. 
$$
(\E^*)^{\otimes k}=\E^*{\otimes}\cdots {\otimes} \E^*
$$
are distributions on the bundle $E \boxtimes \cdots \boxtimes E $ over $X\times \cdots \times X$.  $\Sym^k(\E^*)$ is defined similarly (taking care of the sign of permutations). Then 
$$
   \OO(\E):=\prod_{k\geq 0} \Sym^k(\E^*)
$$
is the analogue of $\OO(V)$.

Since $\omega$ is given by  integration, the Poisson kernel $K_0=\omega^{-1}$ is given by $\delta$-function. In particular, $K_0$ is \emph{singular}. Therefore the naive analogue of the finite dimensional model
$$
  \Delta_{K_0}: \OO(\E)\to \OO(\E)
$$
is \emph{ill-defined} since we can not pair a distribution with another distribution. This difficulty is related to the ultra-violet problem in quantum field theory. Note that the BV bracket $\{-,-\}$ is in fact well-defined for \emph{local functionals}, i.e., functionals of the form $\int_X \mc L$ by integrating a langrangian density. In particular, the classical master equaition $QI_0+{1\over 2}\{I_0, I_0\}=0$ is well-defined for a local interaction $I_0$. 

In \cite{Kevin-book}, Costello uses  a smooth kernel $K_r$ that is homologous to $K_0$ to define the renormalized BV operator (the existence of such a homologous smooth kernel is a consequence of elliptic regularity).  Since smooth functions can be paired with distributions, the BV operator (defined similarly as the finite dimensional case)
$$
  \Delta_{K_r}: \OO(\E)\to \OO(\E)
$$
is well-defined. In this way we obtain a family of differential BV algebras 
$$
\text{BV}[r]:=\bracket{\OO(\E), Q, \Delta_{K_r}}.
$$
Any two such smooth kernels $K_{r_1}$ and $K_{r_2}$ are homologous again by a smooth kernel, which well-defines a homotopic renormalization group flow between $\text{BV}[r_1]$ and $\text{BV}[r_2]$. This motivates the following definition. 

\begin{defn}[\cite{Kevin-book}] An effective solution of quantum master equation is an assignement $I[r]\in \OO(\E)[[\hbar]]$ for each smooth kernel $K_r$ that is homologous to $K_0$, satisfying
\begin{itemize}
\item Renormalized quantum master equation 
$$
(Q+\hbar \Delta_{K_{r}})e^{I[r]/\hbar}=0.
$$
\item Two different $I[r]$'s are related by the homotopic renormalization group flow equation. 
\end{itemize}
If we start with a classical local interaction $I_0$ satisfying classical master equation, then we require that $I[r]$ becomes asymptotically local as $r\to 0$ and  the classical unrenormalized limit $\lim\limits_{r\to 0}\lim\limits_{\hbar\to 0}I[r]=I_0$.  
\end{defn}

There is a parallel observable theory in such a homotopic renormalization framework developed by Costello-Gwilliam \cite{kevin-owen}. Here is a picture to illustrate the situation: 
\begin{center}
  \includegraphics[scale=0.35]{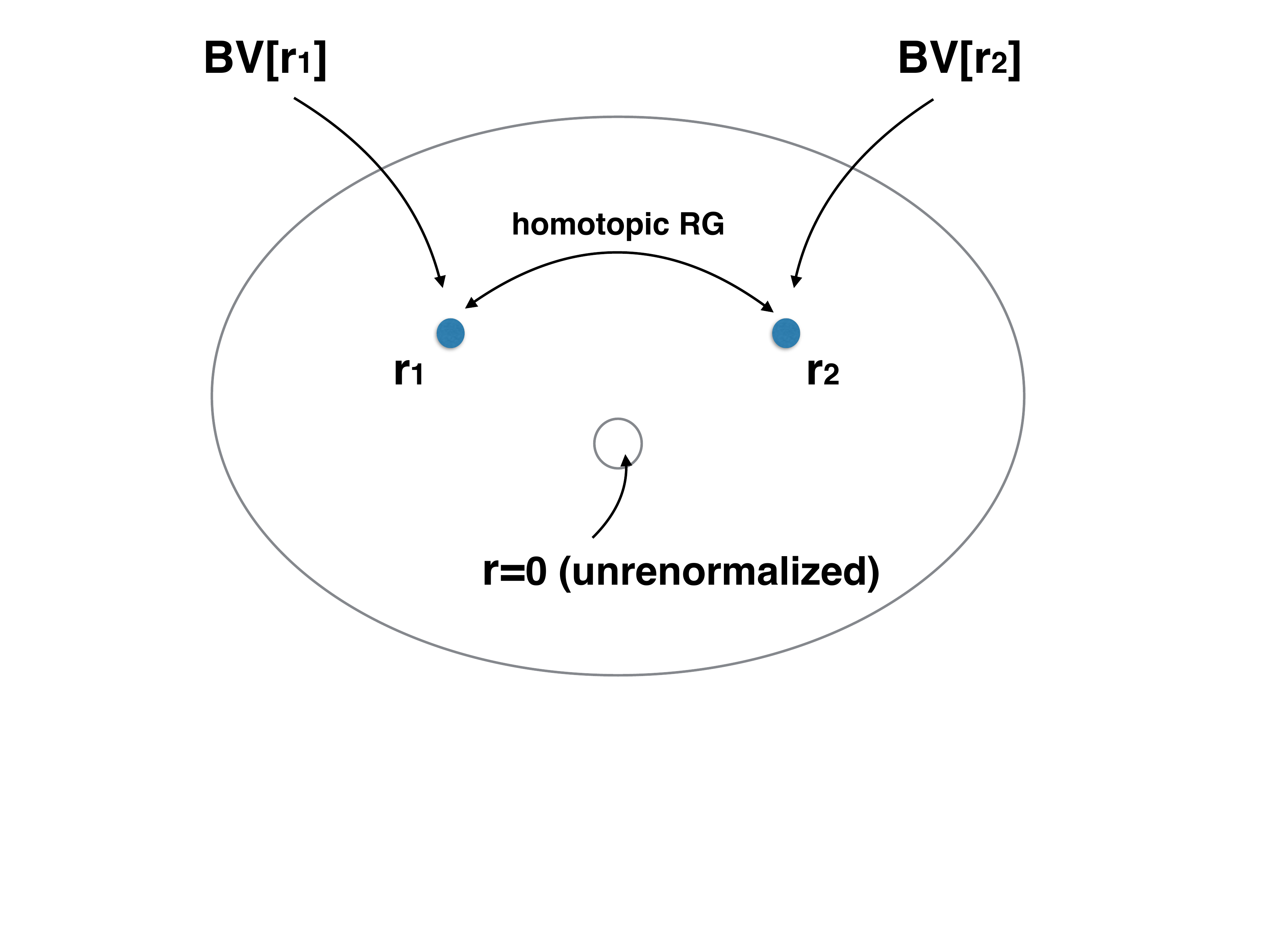}
\end{center}
\begin{rmk}\label{rmk-heat}
Typically for $L>0$,  let $K_L$ be the kernel for the heat operator $e^{-L[Q, Q^{\dagger}]}$ associated to a generlized Laplacian $[Q, Q^{\dagger}]$ constructed from the elliptic operator $Q$ and its adjoint $Q^{\dagger}$. Then $K_0=\lim\limits_{L\to 0}K_L$ is the $\delta$-function which is homologous with $K_L$ by
$$
\int_0^L du (Q^{\dagger}\otimes 1)K_u.
$$
For $ 0<\epsilon<L$, $K_\epsilon$ and $K_L$ are homologous by the smooth kernel
$$
P_\epsilon^L=\int_\epsilon^L du (Q^{\dagger}\otimes 1)K_u
$$
which is called the \emph{regularized propagator}. It leads to the formalism of homotopic renormalization group flow used in \cite{Kevin-book}. 
\end{rmk}
\section{Applications}

\subsection{Topological quantum mechanics and index theorem}
The simplest nontrivial example is when the underlying manifold $X$ is one-dimensional. This corresponds to quantum mechanical models. We give a complete description of the geometry of BV quantization for the low energy effective theory of topological quantum mechanical models. 

Let us start with a linear model. Let $X=S^1$. Let $V$ be a graded vector space  with a degree 0 symplectic pairing
$$
  (-,-): \wedge^2 V\to \R. 
$$ 
The space of fields of our local model will be 
$$
   \E= \Omega^\bullet(S^1)\otimes_{\R} V. 
$$
Let us denote by $X_{dR}$ the super-manifold whose underlying topological space is $X$ and whose structure sheaf is the de Rham complex on $X$. We can identify $\varphi\in \E$ with a map $\hat \varphi$ between super-manifolds
$$
 \hat \varphi: X_{dR} \to V.
$$

The differential $Q=d_{S^1}$ on $\E$ is the de Rham differential on $\Omega^\bullet(S^1)$. The $(-1)$-shifted symplectic pairing is given by the composition of the symplectic pairing and the integration map
$$
  \omega(\varphi_1, \varphi_2):=\int_{S^1} (\varphi_1, \varphi_2), \quad \varphi_1, \varphi_2\in \E. 
$$
This model is just the AKSZ-formalism \cite{AKSZ} applied to the one-dimensional $\sigma$-model. Given $I\in \OO(V)$ of degree $k$, it induces an element $\hat I\in \OO(\E)$ of degree $k-1$ via
$$
   \hat I(\varphi):=\int_{S^1} \hat\varphi^*(I), \quad \forall \varphi \in \E. 
$$

We choose the standard flat metric on $S^1$ and work with the heat kernel regularization as in Remark \ref{rmk-heat}. We have a regularized BV kernel $K_L$ for each $L>0$. $K_\epsilon$ and $K_L$ are homologous by the regularized propagator $P_\epsilon^L$ for $0<\epsilon<L\leq \infty$. The homotopic renormalization group flow operator $W(P_\epsilon^L,-)$ connects the corresponding solutions of quantum master equations. 

\begin{thm}[\cite{GLL}]\label{thm-weyl} Given $I\in \OO^+(V)$ of degree $1$, the limit 
$$
   \hat I[L]= \lim_{\epsilon \to 0} W(P_\epsilon^L, \hat I)
$$
exists as a degree $0$ element of $\OO(\E)[[\hbar]]$. The family $\{\hat I[L]\}_{L>0}$ gives an effective solution of quantum master equation if and only if 
$$
   \bbracket{I, I}_\star=0. 
$$
Here $\OO(V)[[\hbar]]$ inherits a natural Moyal-Weyl product $\star$ from the linear symplectic form $(-,-)$. $[-,-]_\star$ is the commutator with respect to the Moyal-Weyl product. 
\end{thm}

This theorem illustrates an explicit relationship between BV quantization in one dimension and the Weyl quantization. It can be globalized on a symplectic manifold $(X, \omega)$ as follows. For each point $p\in X$, the tangent space $T_pX$ is a linear symplectic space and $\OO(T_pX)$ has a canonical Weyl quantization. It glues to a bundle of algebra $Weyl(T_X)$ with fiberwise Moyal-Weyl product $\star$
$$
  \xymatrix{
     \OO(T_X) \ar[d] \ar@{~>}[rr]^{\text{quantization}}&&  Weyl(T_X)\ar[d]\\
        X &  &X
  }
$$
Let $\nabla$ be a torsion-free connection on $T_X$ that is compatible with $\omega$ (called a symplectic connection).  $\nabla$ naturally induces a connection on $Weyl(T_X)$. In \cite{Fedosov-DQ}, Fedosov shows that $\nabla$ can be modified to be flat by adding a ``quantum correction" $\gamma \in \Omega^1(X, Weyl(T_X))$. Here $\gamma$ satisfies 
$$
    \nabla \gamma+{1\over 2\hbar}[\gamma, \gamma]_\star= \omega_\hbar
$$
where $\omega_\hbar=-\omega+\sum_{k\geq 1}\hbar^k \omega_k$ and $\omega_k's$ are closed 2-forms that classify the equivalence class of deformation quantizations on $X$. Since $\omega_\hbar$ is valued in the central element of $Weyl(T_X)$, the curvature of the quantum corrected connection $\nabla+{1\over\hbar}[\gamma,-]_\star$ vanishes. Fedosov calls this an abelian connection. Let
$$
 \sigma: \Gamma(X, Weyl(T_X))\to \cinfty(X)[[\hbar]]
$$
be the symbol map by projecting to the constant component. Let  $\Gamma^{flat}(X, Weyl(T_X))$ be the flat sections of the Weyl bundle under an abelian connection $\nabla+{1\over\hbar}[\gamma,-]_\star$. Fedosov further shows that 
$$
\sigma: \Gamma^{flat}(X, Weyl(T_X))\to \cinfty(X)[[\hbar]]
$$
is a $\hbar$-linear isomorphism, which induces a star product (deformation quantization) on $\cinfty(X)[[\hbar]]$ via the fiberwise Moyal-Weyl product on $\Gamma^{flat}(X, Weyl(T_X))$.

We consider the topological $\sigma$-model of AKSZ type
$$
   S^1_{dR}\to X. 
$$
We consider the low energy effective theory $\E$ which describes the formal neighborhood of constant maps. It can be viewed as a sheafification of our linear model $V$ by gluing the tangent spaces $\{T_pX\}$. 
\begin{thm}[GLL]\label{thm-Fedosov}
Solutions of effective quantum master equation for $\E$ corresponds to Fedosov's abelian connections. The associative algebra of local quantum observables is identified with Fedosov's deformation quantization. 
\end{thm}

We refer to \cite{GLL} for more precise notions in the above theorem. Physically, this theorem says that the low energy effective theory of topological quantum mechanics is described by Fedosov's abelian connections and deformation quantization. There is a natural lagrangian super submanifold $\mc L$ inside $\E$, which induces a correlation function 
$$
    \int_{\mc L}: \  \Obs^q \to \R((\hbar)). 
$$
Restricting to local quantum observables, it can be identified with the unique normalized trace map on the deformation quantized algebra \cite{Fedosov-DQ, Nest-Tsygan}
$$
   \Tr: \cinfty(X)[[\hbar]]\to \R((\hbar)).
$$
Combining with Atiyah-Bott localization method, we can compute 
$$
 \int_{BV} 1= \int_X e^{-\omega_\hbar/\hbar} \hat A(X).
$$
This is the simplest version of \emph{algebraic index theorem} which was first formulated by Fedosov \cite{Fedbook} and Nest-Tsygan \cite{Nest-Tsygan} as the algebraic analogue of Atiyah-Singer type  index theorem. It leads to the Atiyah-Singer index theorem with the help of the symbol calculus of pseudodifferential operators \cite{G-index, Nest-Tsygan2,Fedesov-PDO}. Our formulation is closely related to  the physics interpretation of Atiyah-Singer index theorem \cite{Alvarez, FW,Witten, Witten-index} via semi-classical analysis of supesymmetric quantum mechanics. When the symplectic manifold $X$ is the total space of a cotangent bundle, there is a closely related result \cite{GG} where the quantization is one-loop exact.

\subsection{2d chiral theory in the BV formalism} We move on to consider quantum field theories in two dimension. For simplicity, we consider a flat surface 
$$
\Sigma=\C, \C^*, \text{or}\ E_\tau=\C/\Z\oplus \Z \tau
$$
where $\tau$ is an element of the upper half plane. 

Here are some examples of free {conformal field theory} (CFT) in two dimension.
\begin{itemize}
\item free boson:  $\int_\Sigma \pa\phi\wedge \bar\pa\phi$. $\quad \phi$ a scalar field.
\item bc-system: $\int_\Sigma b\wedge \bar\pa c$. $\quad b,c$ is a pair of fermions.
\item $\beta \gamma$-system: $\int_\Sigma \beta\wedge \bar\pa \gamma$. $\quad \beta, \gamma$ is a pair of bosons. 
\end{itemize}
We consider effective BV quantization in two dimension for chiral deformation of free CFT of the form:
$$
 S= \text{free CFT}+ I.
$$ 
Here $I$ is a local functional of the form
$$
 I=\int d^2z \mathcal L^{hol}(\pa_z\phi, b, c, \beta, \gamma)
$$
where $\mathcal L^{hol}$ is a lagragian density involving only holomorphic derivatives of $\pa_z\phi, b, c, \beta, \gamma$. 

There is an analogue of Theorem \ref{thm-weyl} describing the effective BV quantization, where the Weyl algebra is replaced by the \emph{vertex algebra}. We give a brief discussion on vertex algebras and refer to \cite{Frenkel,Kac} for details. 

A \emph{vertex algebra} is a vector space $\mc V$ with the structure of  \emph{state-field correspondence } (and other axioms like vacuum, locality, etc.)
\begin{align*}
   \mc V &\to End(\mc V)[[z,z^{-1}]]\\
   A &\to Y(A,z)=\sum_n A_{(n)} z^{-n-1} 
\end{align*}

For simplicity, we will just write $A(z)$ for the operator $Y(A,z)$ associated to $A\in \mc V$. It
defines the \emph{operator product expansion} (OPE)

\begin{minipage}{\linewidth}
      \centering
      \begin{minipage}{0.45\linewidth}
         \begin{equation*}
A(z)B(w)=\sum_{n\in \Z} {(A_{(n)}\cdot B)(w)\over (z-w)^{n+1}}
\end{equation*}
      \end{minipage}
      \hspace{0.05\linewidth}
      \begin{minipage}{0.45\linewidth}

 \includegraphics[scale=0.1]{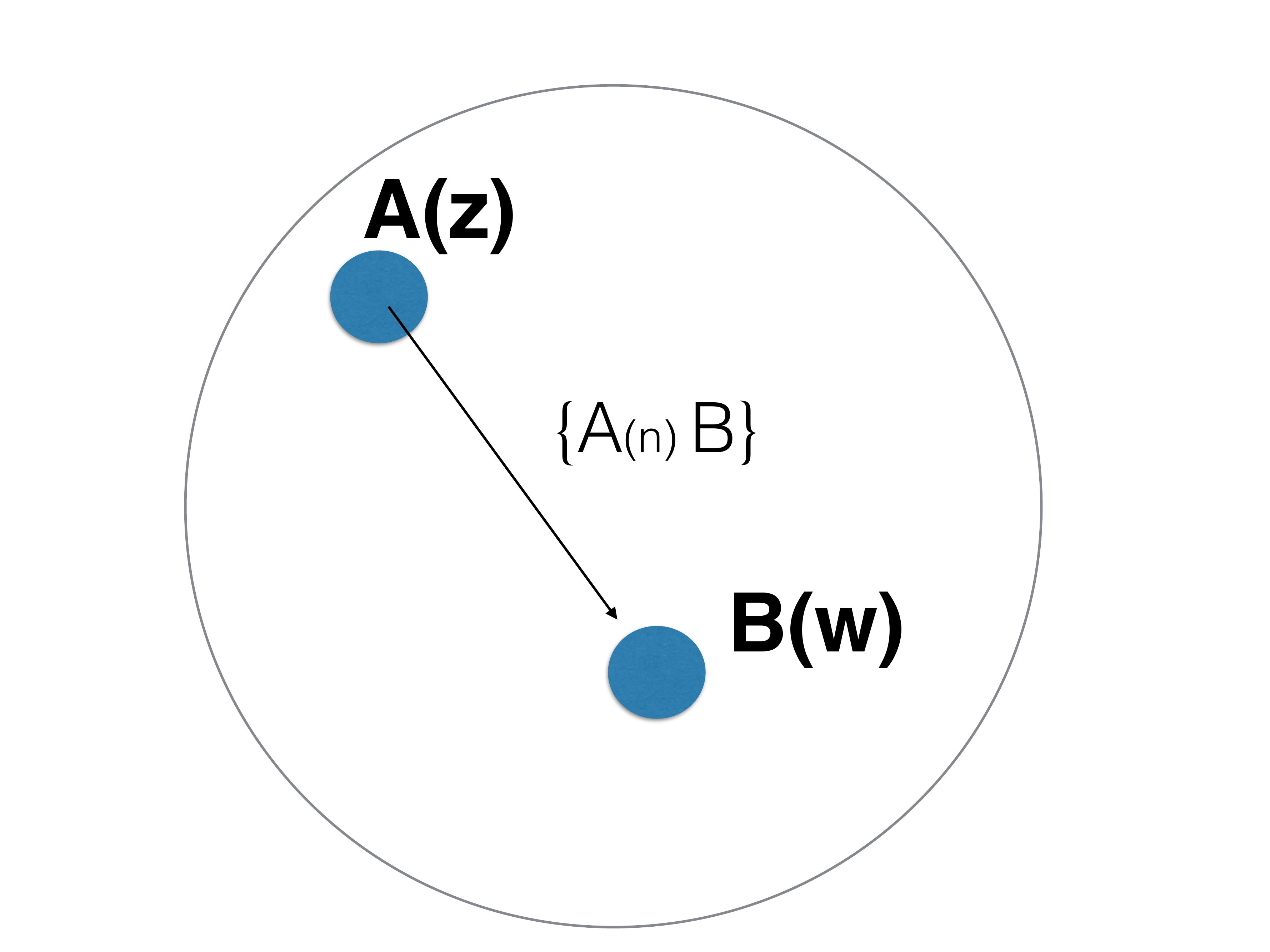}

      \end{minipage}
  \end{minipage}
which is the 2d analogue of ``associative product". We often write its singular part as 
$$
A(z)B(w)\sim \sum_{n\geq 0} {(A_{(n)}\cdot B)(w)\over (z-w)^{n+1}}.
$$

Given a vertex algebra $\mc V$, we can associate a Lie algebra $\oint \mc V$ from its Fourier modes. 
 As a vector space 
$$
  \oint \mc V:=\text{Span}_{\C}\fbracket{\oint dz z^k A(z):=A_{(k)}}_{A\in \mc V, k\in \Z}. 
$$
The Lie bracket is determined by the OPE  (Borcherds commutator formula)
$$
    \bbracket{A_{(m)}, B_{(n)}}=\sum_{j\geq 0}\binom{m}{j}\bracket{A_{(j)}B}_{m+n-j}. 
$$ 
\begin{rmk}\label{rmk-auto}
For any $A\in \mc V$, the residues $A_{(0)}=\oint dz A(z)$ generates an infinitesimal automorphism of the vertex algebra. See for example \cite[Corollary 3.3.8]{Frenkel}.
\end{rmk}

In our BV setup for chiral deformation of free CFT's, the space of fields has the form 
$$
\Omega^{0,\bullet}(E_\tau)\otimes h
$$
of Dolbeault complex valued in some vector space $h$. $h$ generates a vertex algebra $\mc V$, and a chiral lagrangian $\mc L^{hol}$ gives rise to an operator $\oint \mc L^{hol}\in \oint \mc V$.  The following theorem is a 2d analogue of Theorem \ref{thm-weyl}. We refer to \cite{Si-vertex} for details.

\begin{thm}[\cite{Si-vertex}]\label{thm-2d} We consider the following two dimensional quantum field theory with fields $\Omega^{0,\bullet}(E_\tau)\otimes h$ and
$$
  S=\text{free CFT}+ I
$$  
where $I=\int d^2z \mc L^{hol}$ is a {chiral deformation} described above. The differential $Q=\dbar+\delta: \Omega^{0,\bullet}\otimes h\to \Omega^{0,\bullet}\otimes h$ is the sum of $\dbar$ and a holomorphic differential operator $\delta$. 
\begin{enumerate}
\item The theory is {UV-finite}. 
\item Solutions of effective quantum master equations 
$$
\Longleftrightarrow  \delta \oint dz \mc L^{hol}+{1\over 2\hbar}{\bbracket{\oint dz \mc L^{hol}, \oint dz \mc L^{hol}}=0}, 
$$
where $\oint dz \mc L^{hol}$ is viewed as an operator on the corresponding vertex algebra. 
\end{enumerate}
\end{thm}

\begin{rmk} Remark \ref{rmk-auto} implies that we can globalize the above construction to obtain a 2d analogue of Theorem \ref{thm-Fedosov}. For example, chiral de Rham complex \cite{CDR} could be obtained in this way. See also \cite{GGW}.  It would be interesting to obtain a general algebraic index theorem for  chiral vertex operators in this way. 
\end{rmk}

We will explore applications of (2) in later subsections. Here we first explain (1) on the UV finiteness. The UV property for chiral deformations is known to physicists via the method of point-splitting regularization. It has a nice cohomological interpretation in the our current framework. We give an example to illustrate the basic idea. A full development will appear in a forthcoming paper \cite{Si-Jie}. 

Consider the following chiral deformation of free boson on the elliptic curve $E_\tau$
$$
 S={1\over 2}\int_{E_\tau} \pa \phi\wedge \dbar \phi+{1\over 3!}\int_{E_\tau} {d^2z\over \im \tau}  (\pa_z \phi)^3. 
$$
Let us look at a two-loop diagram 

\begin{minipage}{\linewidth}
      \centering
      \begin{minipage}{0.45\linewidth}
         \includegraphics[scale=0.15]{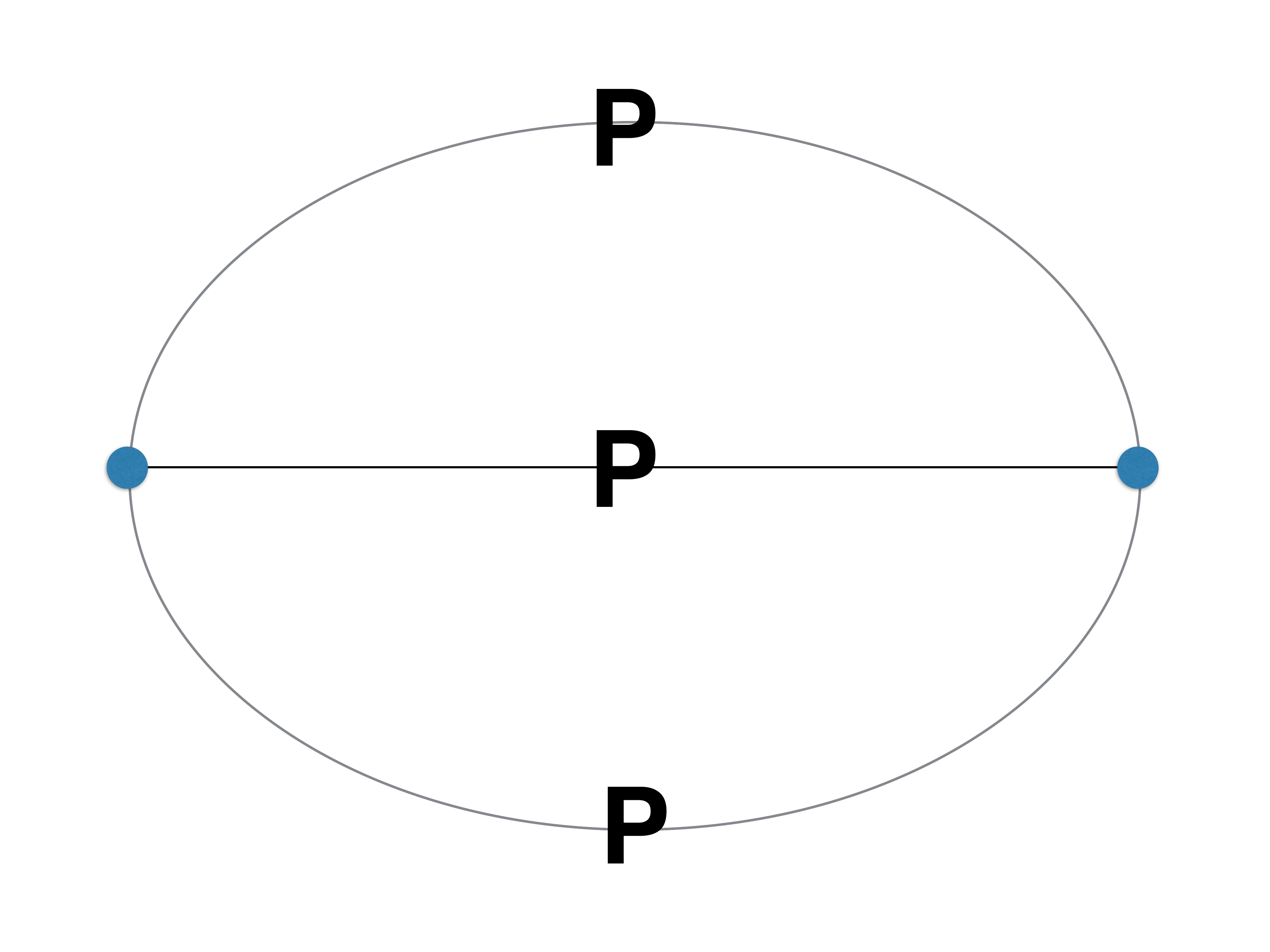} 
      \end{minipage}
      \hspace{-0.2\linewidth}
      \begin{minipage}{0.45\linewidth}
$$
=\int_{E_\tau} {d^2z\over \im \tau}\  \mathbf{P}(z;\tau)^3. 
$$ 
   \end{minipage}
  \end{minipage}
Here the propagator is given by the $\mathbf{P}(z;\tau)=\int_0^\infty {\pa_z^2}h_t$ where $z$ is the linear holomorphic coordinate on $\C$, and $h_t(z)={1\over 4\pi t}\sum\limits_{\lambda\in \Z+\Z \tau}e^{-\abs{z+\lambda}^2/4t}$ is the heat kernel function on $E_\tau$ with respect to the flat metric ${i\over 2}dz\wedge d\bar z$. $\mathbf{P}(z;\tau)$ can be expressed as 
$$
   \mathbf{P}(z;\tau)=\mathcal{P}(z,\tau)+{\pi^2\over 3}E_2^*
$$
where
$\mathcal P$ is the {Weierstrass $\mathcal P$-function}. $E_2$ is the second Eisenstein series, and $E_2^*=E_2-{3\over \pi}{1\over \im \tau}
$.

 Naively, $\mathcal P$ has a pole of order 2, and the above integral looks like divergent. However, our homotopic renormalization scheme leads to a well-defined renormalized value, which can be described as follows. Consider the exact sequence of sheaves on $E_\tau$
 $$
   0\to \C \to  \mathcal M \stackrel{d}{\to} \Omega^1_{II}\to 0
 $$
 where $\mathcal M$ is the sheaf of meromorphic functions, and $\Omega^1_{II}$ is the sheaf of abelian differentials of second kind, i.e. meromorphic 1-forms without order 1 pole. The integrand 
$$
{d^2z\over \im \tau}\  \mathbf{P}(z;\tau)^3=-d\bar z\wedge { \mathbf{P}(z;\tau)^3 dz\over \im \tau}
$$
can be viewed as an element of $H^1(E_\tau, \Omega^1_{II})$. Then the renormalized value of the above Feynman integral is given by the boundary map
$$
  H^1(E_\tau, \Omega^1_{II})\to H^2(E_\tau, \C)\stackrel{\int}{\to} \C. 
$$
Using this we find 
$$
  {1\over \pi^6} \int_{E_\tau} {d^2z\over \im\tau} \mathbf{P}^3={2^2\over 3^3 5}E_6+{2\over 3^2 5}E_4 E_2^*-{2\over 3^3}(E_2^*)^3
$$
Under the $\bar\tau\to \infty$ limit, which amounts to replacing $E_2^*\to E_2$, the above integral reduces to an A-cycle integral computed by Douglas \cite{Douglas}. This explanation of renormalization works in great generality \cite{Si-Jie}. 

\subsection{BCOV theory and higher genus B-model}
The closed string field theory is described in terms of BV master equation \cite{Zwiebach}. This is one of the places where BV quantization plays an essential role. We are interested in mirror symmetry between topological strings. The leading  cubic approximation of string field action in the B-twisted topological sector  on Calabi-Yau 3-fold is proposed in \cite{BCOV}, called the \emph{Kodaira-Spencer gauge theory}. This is  generalized in \cite{Si-Kevin} where  we have found the full B-twisted topological string field action on Calabi-Yaus of arbitary dimension that satisfies the classical master equation. We call this \emph{BCOV theory} and propose (as a generalization of \cite{BCOV}) that the BV quantization of BCOV theory leads to quantum B-model that is mirror to the A-model Gromov-Witten theory of counting higher genus curves. As an application of the effective BV quantization method, we carry out this program explicitly for elliptic curves following \cite{Si-vertex}.  See 
\cite{Partition, TCFT, KS, CT} for another categorical approach to this example. 
 \subsubsection{classical BCOV theory}
Let $X$ be a compact Calabi-Yau manifold of dimension $d$. $\Omega_X$ will be a fixed holomorphic volume form. By abuse of notation, we let
\begin{align*}
\PV(X)=\bigoplus_{0\leq i,j \leq d}\PV^{i,j}(X), \quad
\PV^{i,j}(X)= \Omega^{0,j} (X, \wedge^i T_X)
\end{align*}
be the space of polyvector fields on $X$. Here $T_X$ is the
holomorphic tangent bundle, and $\Omega^{0,j} (X, \wedge^i T_X) $
is the space of smooth $(0,j)$-forms valued in $\wedge^i T_X$. 
 $\Omega_X$ induces an identification between the space of polyvector fields and differential forms
\begin{align*}
   \PV^{i,j}(X) & \stackrel{\lrcorner \Omega_X}{\iso} \Omega^{d-i,j}(X)\\
    \alpha &\to\alpha \lrcorner \Omega_X
\end{align*}
where $\lrcorner$ is the contraction, and $\Omega^{i,j}(X)$ denotes smooth differential forms of type $(i,j)$. The holomorphic de Rham differential $\pa$ on forms defines an operator on $\PV(X)$ via the above isomorphism, denoted again by
$$
\partial : \PV^{i,j}(X) \to \PV^{i-1,j}(X)
$$
i.e.
\begin{align*}
(\partial \alpha) \lrcorner \Omega_X := \partial ( \alpha \lrcorner \Omega_X), \ \ \alpha\in \PV(X).
\end{align*}
The definition of $\pa$ doesn't depend on the choice of $\Omega_X$ on compact Calabi-Yau manifolds.

\begin{defn}[\cite{Si-Kevin}] The space of fields for BCOV theory on $X$ is the complex
$$
  \mc E:= (\PV(X)[[t]][2], Q=\dbar+t\pa)
$$
where $t$ is a formal variable of degree $2$ representing the gravitational descendants. $[2]$ is the degree shifting such that $\PV^{1,1}(X)$ lies in degree  $0$. 
\end{defn}

The BV kernel of BCOV theory is in fact degenerate, hence does not come from a $(-1)$-symplectic pairing. Nevertheless, the framework of effective BV quantization still works (see Remark \ref{rmk-BV-degenerate}). Consider the trace map
$$
  \Tr: \PV(X)\to \C, \quad \mu \to \int_X (\mu\lrcorner \Omega_X)\wedge \Omega_X. 
$$
The BV kernel is a $(\PV(X\times X)=)\PV(X)\otimes \PV(X)$-valued distribution representing the integral kernel (with respect to $\Tr$) of the operator $\pa: \PV(X)\to \PV(X)$. In other words, it is given by 
$$
 (\pa\otimes 1)\delta
$$
where $\delta$ is the $\delta$-function representing the integral kernel of the identity operator with respect to $\Tr$.

 In particular,  $\PV(X)t^k$ for $k>0$ does not appear in the BV kernel:  they are not dynamical  hence are called \emph{background fields}.  We see that BCOV theory is a highly degenerate BV theory with a huge number of background fields. 
 
The classical interaction of BCOV theory is given by  
$$
  I_0(\mu)=\Tr \abracket{e^{\mu}}_0, \quad \mu \in \mc E, 
$$
 where 
 $$
    \abracket{\mu_1 t^{k_1}, \cdots, \mu_n t^{k_n}}_0:=\binom{n-3}{k_1,\cdots,k_n} \mu_1\wedge\cdots \wedge \mu_n, \quad \mu_i\in \PV(X)
 $$
 and $\abracket{e^{\mu}}_0$ is understood as $\sum\limits_{n\geq 3}{1\over n!}\abracket{\mu, \cdots, \mu}_0$. It is shown in \cite{Si-Kevin} that $I_0$ satisfies the classical master equation
 $$
  QI_0+{1\over 2}\{I_0, I_0\}=0. 
 $$
The effective BV quantization of classical BCOV theory leads to higher genus B-model invariants which is expected to be identified with the A-model Gromov-Witten invariants of the mirror Calabi-Yau \cite{BCOV, Si-Kevin}. So far we have only established quantizations of BCOV theory on a few special Calabi-Yau geometries \cite{open-closed, Si-vertex}. Elliptic curve is an example where we have fully established the quantum BCOV theory.

\begin{rmk}The above setup is naturally genearlized to the Landau-Ginzburg case when $X$ is a non-compact Calabi-Yau equipped with a holomorphic function $f: X\to \C$  (superpotential). There $\E=\PV_c(X)[[t]]$ is given by compactly supported polyvector fields, and $Q=\dbar+t\pa+df\lrcorner$. When $X=\C^n$ and $f$ is a holomorphic function with an isolated critical point, the genus zero Landau-Ginzburg twisted BCOV theory is equivalent \cite{LLSaito, LLSS} to K. Saito's theory of primitive forms  \cite{Saito-primitive}. 

\end{rmk} 

\subsubsection{The elliptic curve example}
We carry out the effective BV quantization of BCOV theory explicitly in the case of elliptic curves. 

Let $\tau$ be an element of the upper half plane. Let $\Etau=\C/(\Z\oplus Z\tau)$ be the corresponding elliptic curve. $z$ denotes the linear holomorphic coordinate. The space of fields of BCOV theory on $\Etau$ is 
$$
   \E=\Omega^{0,\bullet}(\Etau, \OO_{\Etau})[[t]]\oplus \Omega^{0,\bullet}(\Etau, T_{\Etau}[1])[[t]]. 
$$
The differential is given by 
$$
  Q=\dbar+t\pa,
$$
where $\pa: T_{\Etau}\to \OO_\Etau$ is the divergence operator with respect to the holormophic volume form $dz$. 

Let us first describe the BV bracket concretely. Let us parametrize a field $\varphi \in \E$ by 
$$
 \varphi=\sum_{k\geq 0}b_kt^k+\eta_k  t^k, \quad b_k \in \Omega^{0,\bullet}(\Etau, \OO_{\Etau}), \quad \eta_k \in \Omega^{0,\bullet}(\Etau, T_{\Etau}[1]).
$$
Then the only nontrivial BV bracket is given by 
$$
   \{b_0(z), b_0(w)\}=\pa_z \delta^2(z-w)
$$
where $\delta^2(z-w)$ is the $\delta$-function  with respect to the volume form ${i\over 2}dz\wedge d \bar z$. All other fields $b_{>0}, \eta_{\bullet}$ are background fields. This theory can be identified with a chiral deformation of the free chiral boson $\phi$, where 
$$
  b_0= \pa_z \phi, \quad \text{with OPE}\quad b_0(z)b_0(w)\sim {\hbar \over (z-w)^2}. 
$$
In particular, Theorem \ref{thm-2d} applies and the effective BV quantization amounts to find an operator $I$ on the vertex algebra satisfying 
$$
   (t\pa) I+{1\over 2\hbar}[I,I]=0
$$
and reduces to our classical BCOV interaction $I_0$ modulo $\hbar$. In \cite{Si-vertex},  a canonical solution of quantum master equation is found using boson-fermion correspondence. The full solution is determined by its restriction to 
$$
     b_{>0}=0, \quad \eta_{\bullet}\in \C[d\bar z]\pa_z\quad  (\text{i.e., $b_{>0}$ vanishes and $\eta_\bullet$=constants}). 
$$
Such subsector represents $t\pa$-cohomology of the fields, and is called the stationary sector \cite{L-elliptic} which is mirror to a similar notion in the A-model. It determines the full generating function on elliptic curves with the help of Virasoro constraint \cite{virasoro}. 

Restricting to the stationary subsector, the quantum corrected chiral interaction is given by
$$
I=\sum_{k\geq 0} \eta_k \int d^2 z {W^{(k+2)}(b_0)\over k+2}
$$
where
$$
  {W^{(k)}(b_0)=\sum_{\sum_{i\geq 1}ik_i=k}{k! \over \prod_i k_i!} :\prod_i\bracket{{ {1\over i!}(\sqrt{\hbar}\pa)^{i-1} b_0}}^{k_i}}:_B
$$
comes from the bosonic realization of the $W_{1+\infty}$-algebra. In the stationary subsector, quantum master equation is equivalent to 
$$
   \bbracket{\oint dz {W^{(k+2)}\over k+2}, \oint dz {W^{(m+2)}\over m+2}}=0, \quad k,m\geq 0,
$$
representing an integrable hierarchy of commuting Hamiltonians. 

The character on the Heisenberg vertex algebra (generated by $b_0$)
$$
\Tr  q^{L_0-{1\over 24}}e^{{1\over \hbar}\sum\limits_{k\geq 0}\oint dz \eta_k {W^{(k+2)}\over k+2}}, \quad q=e^{2\pi i\tau},
$$
coincides with the A-model stationary Gromov-Witten invariants \cite{virasoro} under the boson-fermion correspondence. This can be viewed as a full generalization of \cite{Dijkgraaf-elliptic} who considers the cubic interaction $W^{(3)}$.
 \subsection{Calabi-Yau geometry and integrable hierarchy}
 Topological string on Calabi-Yau manifolds is often found to be related to certain integrable hierarchies. In this section we propose a uniform interpretation in terms of BCOV theory. The classical aspects of this construction is explored in \cite{HSY}.

Let $X$ be a Calabi-Yau geometry. It can be either a compact Calabi-Yau, or a noncompact Calabi-Yau with a Landau-Ginzburg superpotential $f: X\to \C$ having isolated critical points.  We consider $X\times \C$, which is again a Calabi-Yau manifold. Here we equip $\C$ with the trivial volume form $dz$. We consider BCOV theory on $X\times \C$, with a parameter $\epsilon \neq 0$
$$
\mc E^{\epsilon}=\bracket{ \PV(X\times \C)[[t]], Q_\epsilon=\dbar+ t(\pa_X + \epsilon \pa_{\C})}
$$
Here $\pa_X$ ($\pa_{\C}$ respectively) is the divergence operator with respect to the Calabi-Yau volume form on $X$ ($\C$ respectively). The BV kernel is the integral kernel of the operator $(\pa_X + \epsilon \pa_{\C})$ on $\PV(X\times \C)$. It is easy to see that the same classical interaction satisfies the classical master equation. Now we consider the process of integrating out massive modes on $X$
$$
\xymatrix{
  \text{B- model (BCOV theory $\mc E^\epsilon$) on} & X\times \C\ar[d]^{\text{integration out over $X$}}\\
  \text{Effective 2d chiral theory on}& \C
}
$$
Mathematically, this process is a partial homological perturbation along component of $X$ (this requires a choice of splitting of the Hodge filtration on $X$). We end up with an effective 2d chiral theory on $\C$, whose space of fields is
$$
   \PV(\C)\otimes_{\C} \text{H}_X[[t]]
$$
where 
$$
{H}_X=\begin{cases}
 H^\bullet(\PV(X), \dbar) &  \text{compact Calabi-Yau model}\\
 \Jac(f) & \text{Landau-Ginzburg model}.
\end{cases}
$$
$H_X[[t]]$ represents the zero modes of BCOV theory on $X$. The effective theory has again a huge number of background symmetry, which contains the abelian subalgebra 
$$
   \pa_z \otimes \text{H}_X[[t]].
$$
Their associated Noether currents give rise to infinite many commutating operators $\{I_{\alpha}^{\hbar, \epsilon}\}_{\alpha \in \text{H}_X[[t]]}$.  This is our interpretation of integrable hierarchy. 

Observe that the effective quantum parameter on $\C$ is $\hbar \epsilon$. In particular, the limit $\epsilon\to 0$ describes a family of classical master equations parametrized by $\hbar$. 

Let us first consider the limit
$$
 I_{k,\mu}= \lim_{\substack{\hbar\to 0\\ \epsilon\to 0}} I^{\hbar, \epsilon}_{t^k\mu}, \quad \mu \in H_X. 
$$
This can be described as follows. Recall that $H_X$ carries a nondegenerate pairing $\abracket{-,-}$. Our effective vertex algebra is generated by $H_X$ with OPE
$$
   a(z)b(w)\sim {\hbar\epsilon \abracket{a,b}\over (z-w)^2}, \quad a, b\in H_X.
$$
The classical limit leads to the Poisson bracket
$$
  \{a(z), b(w)\}=\abracket{a,b}\pa_z\delta(z-w). 
$$
 Let $F_0^X\in \OO(H_X[[t]])$ be the genus zero generating function on $H_X[[t]]$ of BCOV theory (B-model). Let
$$
     J_{k,\mu}:= \pa_{t^k \mu}F_0^X|_{\H_X}\in \OO(\H_X)
$$
be the generating function with one descendant insertion $t^k\mu$ and all others primary insertions $t^{0}H_X$. Then 
$$
  I_{k,\mu}=\oint dz J_{k,\mu}(b(z)), \quad b\in H_X
$$
and classical master equation implies that they all commute under the Poisson bracket above \cite{HSY}. This is precisely Dubrovin's classical integrable hierarchy associated to a Frobenius manifold \cite{Dubrovin}. The limit
$$
  \lim_{\substack{\epsilon\to 0}} I^{\hbar, \epsilon}_\alpha
$$
gives a dispersion deformation of such a classical integrable hierarchy which is contructed and studied extensively by Dubrovin-Zhang \cite{Dubrovin-Zhang} for semi-simple Frobenius manifold. Topological string on $X\times \C$ suggests that such a dispersion deformation and a further quantization should exist in general. One possible approach is to establish quantum BCOV theory on $X\times \C$. We hope to explore this in the future. 


\address{\tiny YAU MATHEMATICAL SCIENCES CENTER,  TSINGHUA UNIVERSITY, BEIJING, CHINA} \\
\indent \footnotesize{\email{sili@mail.tsinghua.edu.cn}}

\end{document}